# ПОСТАНОВКА ЗАДАЧ КОМБИНАТОРНОЙ ОПТИМИЗАЦИИ ДЛЯ РЕШЕНИЯ ПРОБЛЕМ УПРАВЛЕНИЯ И ПЛАНИРОВАНИЯ ОБЛАЧНОГО ПРОИЗВОДСТВА[1]

М.В. Сарамуд[1], Е.А. Спирин[2], Е.П. Талай[3], Я.Ю. Пикалов[4]

*Сибирский государственный университет науки и технологий имени академика М.Ф. Решетнева, Россия, Красноярск, [1] msaramud@gmail.com, [2] Spirin-Evgeniy@yandex.ru, [3] talay.work.45@gmail.com, [4] yapibest@mail.ru*

*Аннотация.* Рассмотрено применение задач комбинаторной оптимизации к решению проблем планирования процессов для производств, основанных на фонде переналаживаемых производственных ресурсов. Приведены результаты их решения методами смешанного целочисленного программирования.

*Ключевые слова:* оптимизация, планирование на производстве, целочисленное программирование, комбинаторные задачи, облачное производство.

# FORMULATION OF PROBLEMS OF COMBINATORIAL OPTIMIZATION FOR SOLVING PROBLEMS OF MANAGEMENT AND PLANNING OF CLOUD PRODUCTION

M.V. Saramud[1], E.A. Spirin[2], E.P. Talay[3], I.I. Pikalov[4]

*Reshetnev Siberian State University of Science and Technology, Russia, Krasnoyarsk,
[1] msaramud@gmail.com, [2] Spirin-Evgeniy@yandex.ru, [3] talay.work.45@gmail.com, [4] yapibest@mail.ru*

*Abstract.* The application of combinatorial optimization problems to solving the problems of planning processes for industries based on a fund of reconfigurable production resources is considered. The results of their solution by mixed integer programming methods are presented.

**Для цитирования:** Сарамуд М.В., Спирин Е.А., Талай Е.П., Пикалов Я.Ю. Оптимизация организации технологических процессов в условиях облачного производства // Математические методы в технологиях и технике. 2022.

Облачное производство (Cloud-based design and manufacturing, CBDM) - новая производственная парадигма, основанная на фонде переналаживаемых производственных ресурсов, обеспечивающих быстрое получение готового продукта с минимальными затратами на менеджмент. Подобный подход требует быстрой адаптации производства к широкой номенклатуре изделий и построение оптимального плана работ. Гибкости производственных ресурсов может быть обеспечена совокупностью производственных ячеек, содержащих: промышленных роботов; оборудование и оснастку с числовым программным управлением (ЧПУ); приемные приспособления для инструментов; и другое специализированное управление. При этом скорость переналадки производства может быть ограничена скоростью принятия проектных и организационных решений.

Прикладному применению канонических методов решения задач оптимизации производственных процессов, посвящена настоящая работа.

При оптимизации производственных процессов глобально преследуются цель повышения прибыльности предприятия при выпуске продукции, за счет: интенсификации процессов и сокращения времени производственного цикла;

---

снижения издержек при производстве (использование более эффективных конструкторских, технологических, управленческих решений/методов/средств); повышения конкурентоспособности продукции благодаря улучшению характеристик[1].

Описанная цель обычно достигается использованием одного из известных математических инструментов, как правило на некотором временном отрезке, характеризуемом постоянством входных параметров. Облачное производство является динамической системой, которая должна адаптироваться под изменяющиеся внешние и внутренние факторы (набор доступных заказов, поставщиков комплектующих, цен на покупные изделия, состав и количество технологических единиц оборудования и т.д.), поэтому актуальной является задача разработки адаптивной системы управления производственными процессами, использующей комбинация математических инструментов для решения нескольких проблем планирования.

Рассмотрим, как задачи производственного планирования могут быть сведены к решению задачи смешанного целочисленного программирования, среди которых можно выделить следующие комбинаторные задачи:

1. «Задача о ранце» – задача о размещении вещей с различной ценностью, обеспечивающем максимальную ценность выборки;

2. «Задача коммивояжера» – задача поиска выгодного маршрута, обеспечивающего кратчайшей путь;

3. «Задача планирования проекта с ограниченными ресурсами» (Resource Constrained Project Scheduling (RCPSP)) – задача поиска возможного планирования для набора заданий с учетом ограниченности ресурсов и приоритета. Обеспечивает минимальный такт выпуска продукции и максимальное использование доступных ресурсов;

4. «Задача планирования расписания» (Job Shop Scheduling Problem) - упорядочение работ, обеспечивающие минимальный такт выпуска продукции и максимальную загрузку оборудования;

5. «Задача о расположении заводов с нелинейной стоимостью» (Plant Location with Non-Linear Costs);

Также эффективно могут быть применены задачи: «одномерного и двумерного раскроя» (The One-dimensional Cutting Stock Problem, Two-Dimensional Level Packing), решающие проблемы планирования прутковой заготовки для производства деталей различной длины на продольно-токарных автоматах, и листового проката; «распределения частот» (frequency assignment problem), решающая задачу распределения производственных потоков [1]; «о восьми ферзях».

Для решения задачи смешанного целочисленного программирования существует множество как открытых бесплатных реализаций алгоритмов, так и коммерческих решений (GUROBI [2], Octave, Matlab). Рассмотрим алгоритм прикладного применение задачи оптимизации на примерах.

Задача о ранце.

На оборудовании имеется свободное время, которое необходимо использовать с максимальной эффективностью, исходя из имеющейся номенклатуры изделий: количество изделий в номенклатуре $I$ = 6; прибыль от реализации каждого из изделий: $p$ = [10, 13, 18, 31, 7, 15]; продолжительность технологических процессов каждого из изделий: $w$ = [11, 15, 20, 35, 10, 33]; доступное время на оборудовании: $c$ = 47.

Математическая формулировка:

$$\text{maximize}: \sum_{i \in I} p_i \cdot x_i$$
$$\text{subject to}: \sum_{i \in I} w_i \cdot x_i \leq c$$
$$x_i \in \{0,1\} \; \forall i \in I$$

для набора из I сущностей необходимо получить выборку максимальной стоимости, но не превышающую по продолжительности изготовления доступное время на оборудовании.

Полученное решение: наибольшую прибыль *p*=10+31=41 обеспечат изделия 1, 4, с общей продолжительность изготовления *w*=11+35=46. Любая другая комбинация обеспечит снижение прибыли.

Задача коммивояжера.

Логистическая оптимизация транспортных потоков, заключающая в сборе готовых изделий транспортным средством с производственных ячеек: наименования производственных ячеек places = ['СКЛАД','ПЯ1','ПЯ2','ПЯ3','ПЯ4','ПЯ5','ПЯ6','ПЯ7', 'ПЯ8', 'ПЯ9','ПЯ10','ПЯ11','ПЯ12','ПЯ13']; длина путей между ячейками задана треугольной матрицей *dists*.

dists = [[83, 81, 113, 52, 42, 73, 44, 23, 91, 105, 90, 124, 57],
[161, 160, 39, 89, 151, 110, 90, 99, 177, 143, 193, 100],
[90, 125, 82, 13, 57, 71, 123, 38, 72, 59, 82],
[123, 77, 81, 71, 91, 72, 64, 24, 62, 63],
[51, 114, 72, 54, 69, 139, 105, 155, 62],
[70, 25, 22, 52, 90, 56, 105, 16],
[45, 61, 111, 36, 61, 57, 70],
[23, 71, 67, 48, 85, 29],
[74, 89, 69, 107, 36],
[117, 65, 125, 43],
[54, 22, 84],
[60, 44],
[97],
[]]

Математическая формулировка:

$$\text{minimize}: \sum_{i \in I, j \in I} c_{i,j} \cdot x_{i,j}$$
$$\text{subject to}:$$
$$\sum_{j \in V\{i\}} x_{i,j} = 1 \; \forall i \in V$$
$$\sum_{j \in V\{j\}} x_{i,j} = 1 \; \forall j \in V$$
$$y_i - (n+1).x_{i,j} \geq y_i - n \; \forall i \in V\{0\}, j \in V\{0, i\}$$
$$x_{,ji} \in \{0,1\} \; \forall i \in V, j \in V$$
$$y_i \geq 0\{0,1\} \; \forall i \in V$$

для заданного набора точек $V = \{0,..,n-1\}$, расположенных на расстоянии, заданном матрицей *dist*, решение состоит из пар, образующих маршрут движения.

Полученное решение: кратчайший путь длиной 547 обеспечивается следующей траекторией движения транспортного средства: СКЛАД -> ПЯ8 -> ПЯ7 -> ПЯ6 -> ПЯ2 -> ПЯ10 -> ПЯ12 -> ПЯ3 -> ПЯ11 -> ПЯ9 -> ПЯ13 -> ПЯ5 -> ПЯ4 -> ПЯ1 -> СКЛАД.

Задача планирования проекта с ограниченными ресурсами (Resource Constrained Project Scheduling (RCPSP)).

Выполнить планирование взаимосвязанных работ в условия ограниченных ресурсов (но возобновляемых): дана последовательность из 10 взаимосвязанных задач (рис. 1), для каждой из которых определены длительность выполнения и количество требуемых ресурсов, отнесенных к единице времени. Необходимо выполнить планирование задач, обеспечивающее минимальное суммарное время, с ограничением на доступное в каждый момент времени количество ресурсов, соответственно: $c1=6$, $c2=8$.

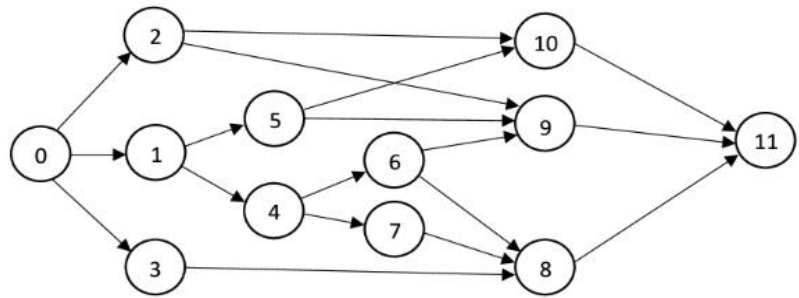

Рис. 1. Последовательность задач, оптимизируемая с учетом ограниченности ресурсов (задачи 0, 11 не требуют ресурсов и символизируют начало и конце процесса)

Таблица 1. Характеристики распределяемых задач

| Работа | Длительность | Требуемое количество ресурса 1 | Требуемое количество ресурса 2 |
|---|---|---|---|
| Job | $p_j$ | $u_{j,1}$ | $u_{j,2}$ |
| 1 | 3 | 5 | 1 |
| 2 | 2 | 0 | 4 |
| 3 | 5 | 1 | 4 |
| 4 | 4 | 1 | 3 |
| 5 | 2 | 3 | 2 |
| 6 | 3 | 3 | 1 |
| 7 | 4 | 2 | 4 |
| 8 | 2 | 4 | 0 |
| 9 | 4 | 5 | 2 |
| 10 | 6 | 2 | 5 |

Математическая формулировка:

$$minimize: \sum_{t \in \tau} t \cdot x_{(n+1,t)}$$
$$subject\ to:$$
$$\sum_{t \in \tau} x_{(i,t)} = 1 \ \forall j \in J$$
$$\sum_{j \in J} \sum_{t_2 = t-p_j+1}^{t} u_{(j,r)} x_{(j,t_2)} \leq c_r \ \forall t \in \tau, r \in R$$
$$\sum_{t \in \tau} t \cdot x_{(s,t)} - \sum_{t \in \tau} t \cdot x_{(j,t)} \geq p_j \ \forall (j,s) \in S$$
$$x_{(j,t)} \in 1 \forall j \in J, \forall t \in \tau$$

для наборов задач $J$, ресурсов $R$, взаимосвязей задач $S$, с горизонтом планирования $\tau$, минимизировать суммарную продолжительность выполнения работ с ограничением $c$ на возобновление ресурсов и заданными продолжительностью работ и потреблением ресурсов

Полученное решение: результатом оптимизации является циклограмма выполнения задач (рис. 2), обеспечивающая минимальное суммарное время выполнения, и максимально допустимое использование возобновляемых ресурсов.

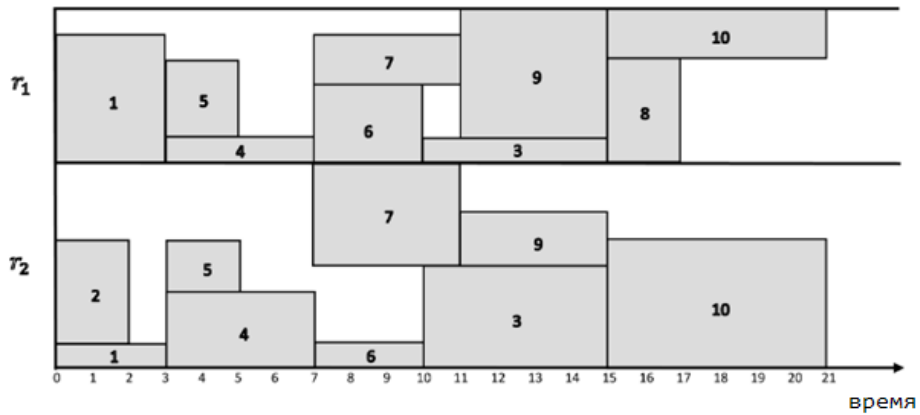

Рис. 2. Циклограмма оптимизированного планирования

Задача планирования расписания оборудования (Job Shop Scheduling Problem).

Необходимо определить последовательность выполнения операций на станках, обеспечив наименьшее время изготовления изделий: количество станков $m=3$; количество изделий $n=3$; длительность выполнения каждой операции на каждом из станков $times$ = [[2, 1, 2], [1, 2, 2], [1, 2, 1]]; последовательность выполнения операций для каждого изделия $machines$ = [[3, 1, 2], [2, 3, 1], [3, 2, 1]].

Математическая формулировка:

$$minimize: C$$
$$subject\ to:$$
$$x_{\sigma_r^j,j} \geq x_{\sigma_{r-1}^j,j} + p_{\sigma_{r-1}^j,j} \forall r \in \{2,..m\}, j \in J$$
$$x_{i,j} \geq x_{i,k} + p_{i,k} - M \cdot y_{i,j,k} \forall j, k \in J, j \neq k, i \in M$$
$$x_{i,k} \geq x_{i,j} + p_{i,j} - M \cdot (1 - y_{i,j,k}) \forall j, k \in J, j \neq k, i \in M$$
$$C \geq x_{\sigma_r^j,j} + p_{\sigma_m^j,j} \forall j \in J$$
$$x_{i,j} \geq 0 \ \forall j \in J, \forall i \in M$$
$$x_{i,j,k} \in \{0,1\}) \forall j, k \in J, i \in M$$
$$C \geq 0$$
$$y_{i,j,k} = \begin{cases} 1, \text{если задача } j \text{ предшествует задаче } k \text{ на станке } i \\ 0, \text{иначе} \end{cases}$$

для множеств задач J и станков M обеспечить минимальную продолжительность выполнения всех операций $C$, с учетом ограничений на производительность $p_{i,j}$ каждого станка $\sigma_m^j$ при выполнении каждой из доступных задач.

Полученное решение: наименьшая длительность изготовления всех изделий равна 7 (рис. 3).

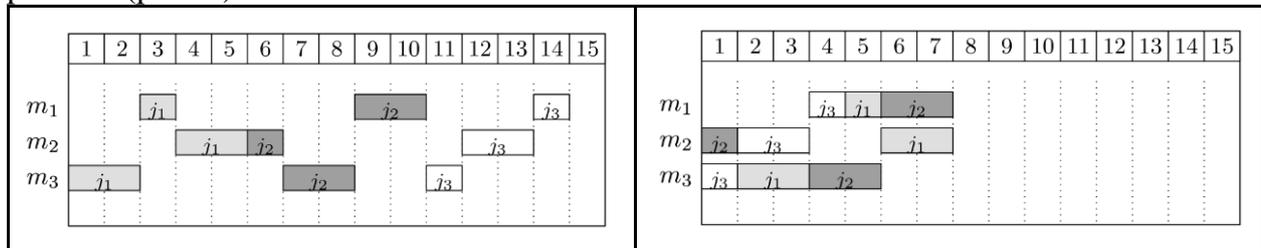

Рис. 3. Наивная (слева) и оптимизированная (справа) циклограммы работы станков

Задача о расположении завода (Plant Location with Non-Linear Costs).

Необходимо расположить два завода в двух регионах (стоимость строительства завода в регионах отличается), удовлетворяющие запросы потребителей, но обеспечивающие минимальную стоимость доставки: идентификаторы возможных мест установки завода $F = [1, 2, 3, 4, 5, 6]$; возможные координаты размещения завода $pf = \{1: (1, 38), 2: (31, 40), 3: (23, 59), 4: (76, 51), 5: (93, 51), 6: (63, 74)\}$; производительность каждого завода в каждом из мест $c = \{1: 1955, 2: 1932, 3: 1987, 4: 1823, 5: 1718, 6: 1742\}$; идентификаторы потребителей $C = [1, 2, 3, 4, 5, 6, 7, 8, 9, 10]$; координаты размещения потребителей $pc = \{1: (94, 10), 2: (57, 26), 3: (74, 44), 4:(27, 51), 5: (78, 30), 6: (23, 30), 7: (20, 72), 8: (3, 27), 9: (5, 39), 10: (51, 1)\}$; объем потребностей потребителей $d = \{1: 302, 2: 273, 3: 275, 4: 266, 5: 287, 6: 296, 7: 297, 8: 310, 9: 302, 10: 309\}$.

Задача решается методом специального упорядоченного набора (SOS [3]) - это упорядоченный набор переменных, используемый как дополнительный способ указать условия целочисленности в модели оптимизации.

Полученное решение представлено на рис. 4:

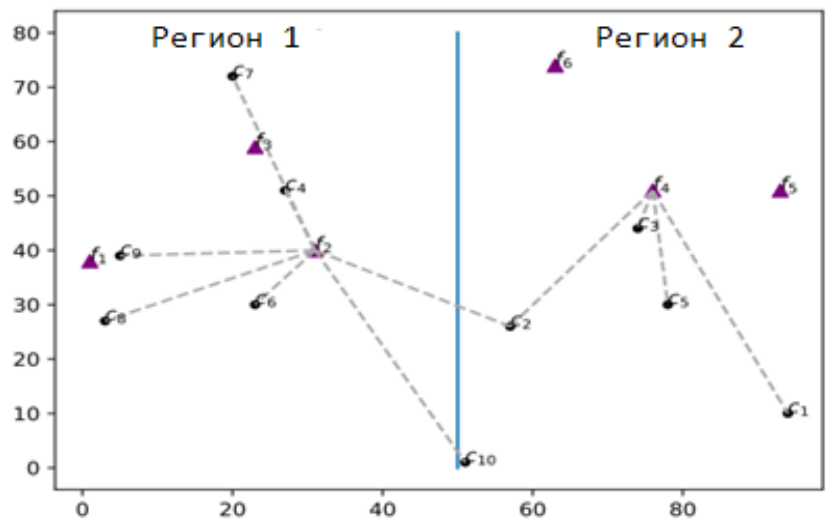

Рис. 4. Наиболее оптимальные из доступных места расположения заводов

Система управления и планирования облачного производства, должна обладать свойствами адаптивности, поэтому может быть построена на комбинации открытых или коммерческих алгоритмов смешанного целочисленного программирования, позволяющих эффективно решать задачи планирования производства.

Библиографический список